\documentstyle{amsppt}
\widestnumber\key{D,G\&G}
%\magnification=\magstep2
%%%makroer
\def\A#1{\Cal A(#1)}%approximabel operators

\def\bb{\beta}
\def\brf{\bar\phi}
\def\bB{\bold B}
\def\bC{\bold C}
\def\bd{\delta}
\def\cB{\Cal B}%Banach algebra of bounded operators
\def\cH{\Cal H}%bounded Hochschild (co)ho
\def\cU{\Cal U}%bounded group of matrices
\def\Diag#1#2{\Delta_{k,N}}
\def\fA{\frak A}%Banach algebra
\def\fB{\frak B}%Banach supalgebra
\def\fC{\frak C}%Banach supalgebra
%Banach algebra generated by group U
\def\fs{\frak s}
\def\ho#1#2{h^N_{#1,#2}}
\def\LIM#1{\operatorname{LIM}_#1}
\def\lp#1{\ell_p(#1)}
\def\nrm#1{\Vert#1\Vert}
\def\ot{\otimes}
\def\pot{\widehat\otimes}
\def\ptens{\pot\cdots\pot}
\def\tD{\widetilde D}

\def\tens{\otimes\cdots\otimes}
\def\nul{\{0\}}

\topmatter
\title
Bounded Hochschild cohomology of Banach algebras with a matrix-like
structure
\endtitle
\rightheadtext{Cohomology of matrix-like Banach algebras}
\author Niels Gr\o nb\ae k
\endauthor
\abstract Let $\fB$ be a unital Banach algebra. A
projection in $\fB$ which is equivalent to the
identitity may give rise to a matrix-like structure on any two-sided ideal $\fA$
in $\fB$. In this set-up we prove a theorem to the effect that the
bounded cohomology $\cH^n(\fA,\fA^*)$ vanishes for all 
$n\geq1$. The hypothesis of this theorem involve (i) strong
H-unitality of $\fA$, (ii) a growth condition on diagonal matrices in $\fA$,
(iii) an extension of $\fA$ in $\fB$ with trivial bounded simplicial
homology. As a corollary we  show that if $X$ is an infinite dimensional
Banach space with the bounded approximation property, $L_1(\mu,\Omega)$
is an infinite dimensional $L_1$-space, and $\fA$ is the Banach algebra of
approximable operators on $L_p(X,\mu,\Omega)\;(1\leq p<\infty)$, then
$\cH^n(\fA,\fA^*)=(0)$ for all $n\geq0$.
\endabstract
\address Department of Mathematics, University of Copenhagen,
Universitetsparken 5, DK-2100 Copenhagen \O, Denmark.\endaddress
\email gronbaek\@math.ku.dk\endemail
\keywords
Bounded Hochschild cohomology, H-unital, simplicially trivial\endkeywords
\endtopmatter
\document
\subheading{0. Introduction} Recall that a Banach algebra $\fA$ is
called weakly amenable if  bounded derivations $D\colon\fA\to\fA^*$
are all inner. Weak amenability of Banach algebras has been the object of
many studies. Recent studies have focused on algebras of compact, in
particular approximable,
operators on Banach spaces. In [D,G\&G] it is shown that if $Y$ is a
reflexive Banach space with
the approximation property, then $\A{\ell_p(Y)},\;1<p<\infty$ is
weakly amenable.  In [B2, Theorem 4.1] it is shown that if $E$
is a  Banach space such that the dual space $E^*$ has the bounded
approximation property and $(\Omega,\Sigma,\mu)$ is a measure space
such that $L_p(\mu)$ is infinite-dimensional, then
$\A{L_p(\mu,E)},\;1\leq p<\infty$ is weakly amenable. In order to prove
this (and many other results on weak amenability) Blanco introduced
so-called trace unbounded systems.  The crucial technical step in
these proofs consists of taking certain averages with respect to such
trace-unbounded systems by exploiting a matrix-like structure. In [G3]
we have shown that $\A{L_p(\mu,E)},\;1\leq p<\infty$ is weakly amenable under
the relaxed assumption (actually both a necessary and sufficient
condition) that $\A{L_p(\mu,E))}$ is {\it self-induced}. 
Loosely speaking this means that  multiplication is an isomorphism as a
balanced bilinear map. In proving this we made use of the matrix-like
structure provided by the right and left shifts on $\ell_p(E)$.

The importance of self-inducedness is that it allows for extending
derivations to enveloping algebras. The definition of self-inducedness
(see [G1]) is equivalent to the part of Wodzicki's definition of
H-unitality [W1]) involving degree 1. Since weak amenability can be
described as the property $\cH^1(\fA,\fA^*)=\{0\}$, it is natural to
investigate the relation between H-unitality and the vanishing of
higher cohomology groups $\cH^n(\fA,\fA^*)=\{0\}$. This is precisely
the concern of the present paper. Specific instances of the question
of vanishing of higher cohomology has been of great importance,
notably the result that $\cH^n(\Cal A,\Cal A^*)=\{0\},\; n\geq0$ for
all $C^*$-algebras without bounded traces [E\&S], and the proof of the
confirmative answer to the 'additive' Karoubi conjecture which
simultaneously gives that $\cH^n(\Cal A,\Cal A^*)=\{0\},\;n\geq0$ for stable
$C^*$-algebras $\Cal A$, see [W2]. For Banach algebras of operators
on general Banach spaces not much is known. The proof in [W2] that
$\Cal B(H)$, the bounded operators on a Hilbert space, has $\cH^n(\Cal
B(H),\Cal B(H)^*)=\{0\},\; n\geq0$ can easily be adapted to replace
$H$ with $\ell_p(X)$, for arbitrary Banach spaces $X$. However, the
method relies on the possibility to make infinite amplifications of
operators, using $H\sim \oplus_1^\infty H$ and thus cannot be adapted
to algebras of compact operators. Our investigations are partly aimed
at circumvening this problem. We prove the following theorem:
 
\proclaim{Theorem} Let $X$ be a Banach space, let $1\leq p<\infty$ and
suppose that the algebra of 
approximable operators $\A{\ell_p(X)}$ is strongly H-unital as a
Banach algebra . Then it is simplicially trivial, i.e.
$\cH^n(\A{\ell_p(X)},\A{\ell_p(X)^*})=\{0\}$ for all $n\geq0$.
\endproclaim

As a corollary we show that all $\A{L_p(\mu,X)},\;1\leq p<\infty$  are
simplicially trivial, when $L_p(\mu)$ is infinite-dimensional and $X$
has the bounded approximation property.

The Banach algebra $\A{\ell_p(X)}$ may be viewed as a topological
tensor product $\A{\ell_p}\bar\otimes\A X$, that is $\A{\ell_p(X)}$ is
reminiscent of a stabilization in $C^*$-algebra theory. It is
therefore natural to expect that, in combination with H-unitality, the
(co)homology is trivial. As Wodzicki's result uses specific
$C^*$-properties, we cannot appeal directly to his approach. In stead
we use that a dense subalgebra has trivial algebraic homology and extend this
by continuity, using elementary chain homotopy techniques.

\subheading{1. The set-up}
\comment
notation hochschild complex, bar complex, approximable operators
upper left hand corners, h-unitality (strong) and consequences,
splitting of multiplication,  
\endcomment

We start by describing the Hochschild homological setting in which we
shall be working. We  appeal to preliminary knowledge and shall be
rather brief.

Let $A$ be an algebra (not necessarily unital) over
$\Bbb C$. We define {\it $n$-chains $C_n(A)$} and {\it face maps}
$d_i\colon C_n(A)\to C_{n-1}(A),\;i=0,\dots, n$ as 
$$  
C_n(A)=A^{\ot(n+1)}
$$
and
$$
d_i(a_0\tens a_n)=\cases a_0\tens
a_ia_{i+1}\tens a_n,&i=0,\dots,n-1\\
a_na_0\tens a_{n-1},&i=n.\endcases
$$
We consider the {\it bar-complex}
$$
B_{\bullet}(A)\colon\quad 0@<<< A@<\bb<< A\ot A@<\bb<<\dots @<\bb<<
A^{\ot n}@<\bb<<A^{\ot(n+1)}\dots, 
$$
with the {\it bar-differential} on $n$-chains as
$$
\bb=\sum_{i=0}^{n-1}(-1)^id_i.
$$
Likewise we have the {\it Hochschild complex}
$$
C_{\bullet}(A)\colon\quad 0@<<< A@<\bd<< A\ot A@<\bd<<\dots @<\bd<<
A^{\ot n}@<\bd<<A^{\ot(n+1)}\dots, 
$$
with the {\it Hochschild differential} on $n$-chains as
$$
\bd=\sum_{i=0}^n(-1)^id_i=\bb+(-1)^nd_n.
$$
For a Banach algebra we consider the corresponding {\it bounded} complexes based on the
projective tensor product.
$$
\bB_{\bullet}(\fA)\colon\quad 0@<<< \fA@<\bb<< \fA\pot \fA@<\bb<<\dots @<\bb<<
\fA^{\pot n}@<\bb<<\fA^{\pot(n+1)}\dots, 
$$
and
$$
\bC_{\bullet}(\fA)\colon\quad 0@<<< \fA@<\bd<< \fA\pot \fA@<\bd<<\dots @<\bd<<
\fA^{\pot n}@<\bd<<\fA^{\pot(n+1)}\dots, 
$$
with $n$-chains $\bC_n(\fA)=\fA^{\pot(n+1)}$ and the same defining
formulas for face maps $d_i$ and differentials $\bb$ and $\bd$ . Upon taking
homology of the complexes $C_{\bullet}(A)$ and $\bC_{\bullet}(\fA)$ we arrive at the
{\it Hochschild homology groups} $H_n(A,A),\;n=0,1,\dots$ respectively
the {\it bounded Hochschild homology groups} $\cH_n(\fA,\fA),\;n=0,1,\dots$.

The duals of the complexes give co-complexes and corresponding
cohomology. Specifically we have 
$$
C^{\bullet}(A)\colon\quad 0@>>> A'@>\partial>> (A\ot
A)'@>\partial>>\dots @>\partial>> (A^{\ot n})'@>\partial>>(A^{\ot(n+1)})'\dots,
$$
and
$$
\bC^{\bullet}(\fA)\colon \quad0@>>> \fA^*@>\partial>> (\fA\pot
\fA)^*@>\partial>>\dots @>\partial>> (\fA^{\pot
  n})^*@>\partial>>(\fA^{\pot (n+1)})^*\dots, 
$$
where $(\cdot)'$ denotes algebraic dual and $(\cdot)^*$ denotes
topological dual. We have used $\partial$ as the common symbol for
algebraic and topological dual maps of the likewise commonly denoted
Hochschild differentials. Correspondingly we have the {\it Hochschild
  cohomology groups} $H^n(A,A')$ and the {\it bounded Hochschild
  cohomology groups} $\cH^n(\fA,\fA^*)$. 
 
We now define the
central notion of H-unitality due to M. Wodzicki [W1]. 

\definition{1.1 Definitions} A complex algebra $A$ is {\it H-unital} if
the bar-complex $B_{\bullet}(A)$ is acyclic. A Banach algebra $\fA$ is
{\it strongly H-unital} if $\bB_{\bullet}(\fA)$ is pure exact,
i.e\. $X\pot\bB_{\bullet}(\fA)$ is acyclic for all Banach spaces
$X$. The Banach algebra $\fA$ is {\it simplicially trivial} if
$\cH^n(\fA,\fA^*)=\{0\}$ for all $n\geq0$.
\enddefinition

The importance of H-unitality is that it is the condition for excission
in Hoch\-schild (co)homology.  Our strategy shall be to use excission
according to a scheme
$$
\CD
0@>>>\fA@>>>\fA\oplus\fC@>>>\fC@>>>0\\
&&@AAA@AAA@AAA\\
0@>>>A@>>>A\oplus C@>>>C@>>>0.
\endCD
$$
where we have control of the algebraic cohomology in the bottom row and
use continuity to extend it to the top row. This is done by
exploiting matrix stuctures (taylored to deal with Banach algebras of
operators on Banach spaces of vector valued sequences) which we now discuss.

 Assume that the Banach  algebra $\fA$ has bounded multipliers $S$ and
$R$ such that $RS=\bold1_\fA$ and $SR\neq\bold1_\fA$. Let $\fA_k=\{a\in\fA\mid
R^ka=aS^k=0\}$. Then $\fA_1\subseteq\fA_2\subseteq\dots$.
Setting $P_i=S^{i-1}R^{i-1}-S^iR^i,\;i=1,2,\dots$ each element
 in $\fA_k$ may be represented
as an infinite upper left corner matrix
$$
(a_{ij})=
\pmatrix
a_{11}&\dots&a_{1k}&0&\dots\\
\vdots&\ddots&\vdots&\vdots&\dots\\
a_{k1}&\dots&a_{kk}&0&\dots\\
0&\dots&0&0&\dots\\
\vdots&&\vdots
\endpmatrix
$$
with $a_{ij}=P_iaP_j$ for some $a\in\fA$. In this picture
multiplication in $\fA_k$ is just matrix multiplication.

\definition{Condition M1} If the algebra of 'upper left hand corners' $\fA_\infty$, given by
$$
\fA_\infty=\bigcup_1^\infty\fA_k\;,
$$
is dense in $\fA$, and if 'block diagonal' matrices are bounded
$$
\sup\{\nrm{\sum_{l=1}^nS^{lk}aR^{lk}}\mid k,n\in\Bbb n\}<\infty,\quad a\in\fA,
$$ 
then $\fA$ satisfies condition M1.
\enddefinition

In order to lift cohomology from a dense subalgebra we shall further
be interested in

\definition{Condition M2} There is a dense subalgebra $A\subseteq\fA$
 such that with $A_k=\{a\in A\mid R^ka=aS^k=0\}=A\cap\frak A_k$, the subalgebra
$$
A_\infty=\bigcup_1^\infty A_k\;,
$$
is H-unital and $H^n(A_\infty,A_\infty')=\{0\}\;,n\geq1$.
\enddefinition

 From the matrix picture of
$\fA_\infty$ it follows that scalar matrices may operate on
$\fA_\infty$ by means of matrix multiplication.
The purpose in our setting of requiring H-unitality is that excission enables
us to extend cocycles to certain operations of permutations. This
leads to the requirement

\definition{Condition M3} Let $\frak S(\Bbb N)$ be the full
permutation group of $\Bbb N$ and consider the subgroup 
$$
\frak G=\{\frak p\in\frak S(\Bbb N)\mid \frak p(i)\neq i\text{ for at most
finitely many }i\in\Bbb N\}
$$
of finite permutations of $\Bbb N$. Let $\pmb D$ be the group of
infinite diagonal scalar matrices with $-1$ in at most
finitely many diagonal entries and $+1$ in the remaining diagonal
entries and let $\Cal U$ be the group of
infinite matrices of the form $U=D \Pi_\frak p$, where $D\in\pmb D$
and 
$\Pi_\frak p$ is the permutation matrix corresponding to some $\frak
p\in \frak G$. If there is $C>0$ such that
$$
\sup\{\nrm{(a_{ij})U},\nrm{U(a_{ij})}\mid
U\in \Cal U\}\leq C\nrm{(a_{ij})},\quad (a_{ij})\in\fA_\infty, $$
that is, $\Cal U$ operates boundedly on $\fA_\infty$ by matrix multiplication of the
corresponding  matrices, then $\fA$ satisfies condition M3.

\medskip
Note that for each $\frak p\in \frak G$ and each diagonal matrix
$D\in\pmb D$ there is $D'\in\pmb D$ such
that $D\Pi_\frak p=\Pi_\frak pD'$
\enddefinition

\subheading{2. The theorem}

\proclaim{2.1 Theorem} Let $\frak A$ be a strongly H-unital Banach
algebra satisfying Conditions {\rm M1, M2,} and {\rm M3}. Then
$$
\cH^n(\fA,\fA^*)=\{0\},\quad n\geq1.
$$
\endproclaim

The proof of this will occupy the rest of the two next section. Thus we shall
assume throughout that $\fA$ satisfies the hypotheses. First we make
some reductions of the problem.

Let $D\in\bold C^n(\fA)$ be an $n$-cocycle, i\.e\. $D$ is a bounded
linear functional on $\fA\ptens\fA$ (equivalently, a bounded
$(n+1)$-linear functional on $\fA$) such that $\partial D=0$. By the
first part of condition M2 there is a (not necessarily bounded) $n$-linear
functional $\phi$ such that 
$$
D(a_0\tens a_n)=\phi(\bd(a_0\tens a_n)),\quad a_i\in A,\;i=0,\dots,n.\tag1
$$

We may prove that $D$ is a (bounded) coboundary by proving that $\phi$
can be chosen to be bounded
on $\bd(A\tens A)$. First we exploit H-unitality. By Condition M3 the
group $\Cal U$ of infinite matrices operates boundedly on $A_\infty$
which by M2 and M1 is dense in $\fA$. Thus we may regard $\Cal U$ as a
bounded group of invertible operators on $\fA$. Let $\fC$ be the
Banach subalgebra of $\cB(\fA)$ generated by $\cU$ and
consider the extension
$$
0@>>>\fA@>\iota>>\fA\oplus\fC@>q>>\fC@>>>0.\tag2
$$

\proclaim{2.2 Proposition} Assume that $\fA$ is strongly H-unital. Then
every bounded $n$-cocycle may be extended to a bounded  $n$-cocycle on
the semi-direct product $\fA\oplus\fC$.
\endproclaim
\demo{Proof} Set $\fB=\fA\oplus\fC$. From [] we get a long
exact sequence  
$$
\dots@>>>\cH_n(\fA,\fA)@>\iota_*>>\cH_n(\fB,\fB)@>>>\cH_n(\fC,\fC)
     @>>>\cH_{n-1}(\fA,\fA)@>\iota_*>>\cdots.\tag3
$$
As in the proof of  [J, Proposition 6.1] we see that 
$\cH_n(\fC,\fC)=\{0\},\;n\geq1$ so that the inclusion 
$\fA@>\iota >>\fB$ induces an isomorphism in bounded homology. Hence
the restriction map induces an isomorphism in bounded cohomology,
i.e. for $n\geq1$ every bounded $n$-cocycle on $\fA$ is the restriction of a
bounded $n$-cocycle on $\fB$.
\enddemo

We shall also need (a consequence of) the discrete version of this
theorem. We let $M_\infty$ denote the algebra of infinite matrices
over $\Bbb C$ consisting of 'upper left-hand corners',
$M_\infty=\{(\lambda_{ij})\mid \exists k\in\Bbb N\colon (i\geq
k\text{ or } j\geq k)\implies\lambda_{ij}=0\}$. We further define $\widetilde
M_\infty=M_\infty\oplus\frak D$, where $\frak D$ is the algebra of
infinite diagonal matrices with constant diagonal, so that we have an
augmentation
$$
0@>>>M_\infty@>>>\widetilde M_\infty@>>>\Bbb C@>>>0\,.\tag4
$$
   As above we may regard $\widetilde M_\infty$ as
multipliers on $\fA$. In fact $\widetilde M_\infty$ is a subalgebra of
$\fC$. We set $B=A\oplus \widetilde M_\infty$. Then 

\proclaim{2.3 Proposition} $H_n(B,B)=H^n(B,B')=\{0\}$ for all $n\geq1$.
\endproclaim
\demo{Proof} Consider the extension
$$
0@>>> A@>>>B@>>>\widetilde M_\infty @>>>0.\tag5
$$
By M2 $A$ is H-unital so that we have excission in
discrete homology, that is, we have a long exact sequence
$$
\dots @>>>H_n(A,A)@>>>H_n(B,B)@>>>H_n(\widetilde M_\infty,\widetilde
M_\infty)@>>>H_{n-1}(A,A)@>>>\dots.\tag6 
$$
The tensor trace is a quasi-isomorphism [L, Theorem 1.2.4], so we get
$H_n(M_\infty,M_\infty)\mathbreak =H_n(\Bbb C,\Bbb C)=\{0\}$.  Having local units
$M_\infty$ is H-unital, yielding
$H_n(\widetilde M_\infty,\widetilde M_\infty)\ =\{0\},\;n\geq0$. By the
assumption M2, $H_n(A,A)=\{0\},\;n\geq1$, so excission gives the
desired result for homology and thereby for cohomology.
\enddemo

We return to the bounded $n$-cocycle $D$. We collect above findings in

\proclaim{2.4 Proposition} Assume that $\fA$ is strongly H-unital, and let
$D$ be a bounded $n$-cocycle on $\fA$. Then there is 
a bounded $n$-cocycle $\tD$ on $\fB=\fA\oplus\fC$ extending $D$ and
there is a (not necessarily bounded) $n$-linear functional
$\widetilde\phi$ on $B=A\oplus\widetilde M_\infty$, such that
$$  
\tD(b_0\tens b_n)=\widetilde\phi(\bd(b_0\tens b_n)),\;b_o,\dots,b_n\in
B.
$$
\endproclaim

\bigskip
\subheading{3. The proof} In the following we drop the\ $\widetilde{}\;$'s,
and hence consider $D$ to be a bounded $n$-cocycle on $\fB$ with
corresponding $n$-linear functional $\phi$ on $B$. 
For each $k\in\Bbb N$ we define a bounded grade zero chain map
$\sigma_k\colon \bC_{\bullet}(\fA)\to \bC_{\bullet}(\fA)$ by 
$$  
\sigma_k(a_0\tens a_n)=S^ka_0R^k\tens S^ka_nR^k.\tag7
$$ 
The first step in proving that $D$ is a bounded co-boundary on $\fA$
is 

\proclaim{3.1 Proposition} Let $\tau$ be an {\rm($n-$1)}-cycle in
$C_{n-1}(A_\infty)$. Then for $k\in\Bbb N$ sufficiently large
$$
\vert \phi(\sigma_k(\tau))-\phi(\tau)\vert\leq K \Vert
D\Vert\Vert\tau\Vert,
$$
where $K>0$ is a constant depending only on $n$.
\endproclaim
\demo{Proof} Choose $k\in\Bbb N$ so large that $\tau\in A_k^{\ot
n}$. Let $U\in\fC$ be the unitary order 2 operator given by
the permutation $(1\,\;k+1)\circ(2\,\;k+2)\circ\dots\circ(k\;2k)$ of
coordinates. Then 
$$
S^kaR^k=UaU^{-1},\;a\in A_k.\tag8
$$ 
Since $U$ is of order 2 we may write $U=e^h$ with
$h=\frac{i\pi}2U$. We shall use a typical homotopy argument to finish. Define
a bounded chain map   
$\alpha_t\colon \bC_{\bullet}(\fA)\to \bC_{\bullet}(\fA)$ by
$$
\alpha_t(a_0\tens a_n)=e^{th}a_0e^{-th}\tens e^{th}a_ne^{-th}.\tag9
$$ 
Then $\alpha_0(\tau)=\tau$ and $\alpha_1(\tau)=\sigma_k(\tau)$. We
further have 
$$
\frac d{dt}\alpha_t(\tau)=\alpha_t(\Gamma_h(\tau)),\tag10
$$
where
$$
\Gamma_h(a_0\tens a_n)=\sum_{i=0}^n\dots\ot[ha_i-a_ih]\ot\dots.\tag11
$$
It is well-known [L, Proposition 1.3.3] that $\Gamma_h$ is 0-homotopic. A bounded 
contracting chain homotopy $s\colon \bC_{\bullet}(\fB)\to
\bC_{\bullet+1}(\fB)$ is given by
$$
s(a_0\tens a_n)=\sum_{k=0}^n(-1)^{k+1}a_0\tens a_k\ot h\ot a_{k+1}\tens
a_n.\tag12
$$
Hence
$$
\aligned
\sigma_k(\tau)-\tau&=\int_0^1\frac d{dt}\alpha_t(\tau)\,dt\\
&=\int_0^1\alpha_t(\Gamma_h(\tau))\,dt\\
&=\int_0^1\alpha_t((s\bd+\bd s)\tau)\,dt\\
&=\bd\int_0^1\alpha_t(s(\tau))\,dt
\endaligned\tag13
$$
where the last step follows from $\bd\tau=0$ and elementary properties
of the Bochner integral. Note that $\int_0^1\alpha_t(\tau)\,dt\in
C_{n-1}(A_\infty)$ and that $s(\tau)\in C_n(B)$. On applying $\phi$ we get
$$
\aligned
\vert\phi(\sigma_k(\tau))-\phi(\tau)\vert&=\vert\phi(\bd\int_0^1\alpha_t(s(\tau))\,dt)\vert\\
&=\Vert D(\int_0^1\alpha_t(s(\tau))\,dt)\Vert\\
&\leq n\frac\pi2 e^{(n+1)\pi\Vert U\Vert}\Vert U\Vert\Vert D\Vert\Vert\tau\Vert\,.
\endaligned\tag14
$$
Since the group $\Cal U$ is bounded, we arrive at the conclusion of the
theorem
\enddemo

In order to utilize this we shall reduce the
problem further. We shall several times take limits by means of the
following functional.

\definition{3.2 Definition} Let
$\operatorname{LIM}\colon\ell_\infty\to\Bbb C$ be a norm one Hahn-Banach extension
of the functional $(x_k)\mapsto\lim_{k\to\infty}x_k$ defined on
the closed subspace of $\ell_\infty$ consisting of convergent
sequences. When we need to refer to the indexing of a sequence, we shall
write
$\LIM kx_k$ to denote the value of the functional LIM on the sequence
$(x_k)_{k\in\Bbb N}$.
\enddefinition

Our strategy hinges on

\proclaim{3.3 Lemma} Let $D\colon\fA^{\pot(n+1)}\to\Bbb C$ be a
bounded $n$-cocycle. Suppose that there is a constant $\kappa>0$ and
for each $k\in\Bbb N$ a bounded $n$-linear functional $\psi_k$ such
that $\Vert\psi_k\Vert\leq\kappa$ and
$$
D(a_0\tens a_n)=\psi_k(\bd(a_0\tens a_n)),\;a_0,\dots,a_n\in A_k.
$$
Then $D$ is a bounded coboundary.
\endproclaim
\demo{Proof} Let $a_0,\dots, a_{n-1}\in A_\infty$. For
 $k\in\Bbb N$
sufficiently large we have $a_0,\dots a_{n-1}\in A_k$. We may consequently
define a bounded $n$-linear functional on $A_\infty$ by
$$
\psi(a_0\tens a_{n-1})=\LIM k(\psi_k(a_0\tens
a_{n-1})),\;a_0,\dots,a_{n-1}\in A_\infty.\tag15 
$$
(Note that, even though $\psi_k(a_0\tens a_{n-1})$ may be undefined for
the first terms, $\psi$ is well defined, because
$\LIM k (x_k)=0$ when $(x_k)$ is eventually 0.) With this definition we have
$$
D(a_0\tens a_n)=\LIM k\psi_k(\bd(a_0\tens a_n))=\psi(\bd(a_0\tens a_n))\tag16
$$
for all $a_0,\dots, a_n\in A_\infty$. Since $\Vert\psi\Vert\leq\kappa$ and
since $A_\infty$ is dense in $\fA$, it follows that $D$ is a bounded coboundary.
\enddemo
We are now ready for
\demo{Proof of Theorem 2.1} In order to find the bounded functionals
$\psi_k$ of Lemma 3.3 we shall make use of
certain chain homotopies. For the rest of the proof we fix $k\in\Bbb N$. We first define certain
graded maps. For each $N\in\Bbb N$ we define the grade 0
maps $\Delta_N\colon \bC_\bullet(\fB)\to\bC_\bullet(\fB)$ by
$$
\Delta_N(a_0\tens a_n)=\sum_{l=0}^{N-1}S^{lk}a_0R^{lk}\tens\sum_{l=0}^{N-1}S^{lk}a_nR^{lk},\tag17
$$ 
and for each $l=0,\dots, N-1$ the grade 1 maps $\ho li\colon\bC_\bullet(\fB)\to\bC_{\bullet+1}(\fB)$ by
$$
\ho li(a_0\tens a_n)=\sigma_{lk}(a_0\tens
a_i)\ot P_l\ot\Delta_N(a_{i+1}\tens a_n),\;i=0,\dots,n,\tag18
$$
where $P_l=S^{(l-1)k}R^{(l-1)k}-S^{lk}R^{lk}$ is the projection onto
the $l$'th coordinate block of size $k$.

\remark{Remark} The grade zero maps $\Delta_N$ can be described by means of
the previously defined chain maps $\sigma_l$:
$$
\Delta_N=\sum_{l=0}^{N-1}\sigma_{lk}\tens\sum_{l=0}^{N-1}\sigma_{lk}.\tag19
$$
\endremark

We now show that the hypothesis of Lemma 3.3 can be fulfilled. The
diagonal matrix picture below will be convenient.

 For $a_0,\dots,a_n\in
A_k$ we have
$$
\Delta_N(a_0\tens a_n)=
\pmatrix a_0&&&&\\
&\ddots&&&\\
&&a_0&&\\
&&&0&\\
&&&&\ddots\endpmatrix\tens
\pmatrix a_n&&&&\\
&\ddots&&&\\
&&a_n&&\\
&&&0&\\
&&&&\ddots\endpmatrix\tag20
$$
as upper left hand $N\times N$ corners of $k\times k$ blocks. We observe
that $\Delta_N$ is a chain map when restricted to $A_k\tens A_k$:
$$
\Delta_N(\bd (a_0\tens a_n))=
\bd(\Delta_N(a_0\tens a_n)),\;a_0,\dots,a_n\in A_k.\tag21
$$
Furthermore
$$
\multline
\ho li(a_0\tens a_n)=\\
\pmatrix
0&&&&\\
&\ddots&&&\\
&&a_0&&\\
&&&0&\\
&&&&\ddots\\
\endpmatrix\tens
\pmatrix
0&&&&\\
&\ddots&&&\\
&&\bold1&&\\
&&&0&\\
&&&&\ddots\\
\endpmatrix\ot\\
\pmatrix a_{i+1}&&&&\\
&\ddots&&&\\
&&a_{i+1}&&\\
&&&0&\\
&&&&\ddots\endpmatrix\tens
\pmatrix a_n&&&&\\
&\ddots&&&\\
&&a_n&&\\
&&&0&\\
&&&&\ddots\endpmatrix,
\endmultline\tag22
$$
as upper left hand $N\times N$ corners with $a_0,\dots,a_i$ and
$\bold1$ placed in the $l$'th diagonal $k\times k$ block entries. 

For $a_0,\dots,a_n\in A_k$ we have the relations
$$
d_i\ho lj(a_0\tens a_n)=\cases \ho l{j-1} d_i(a_0\tens a_n),&i<j\\
                  d_i\ho l{i-1}(a_0\tens a_n),&i=j,j+1\\
                  \ho lj d_{i-1}(a_0\tens a_n),&i>j+1\,.
           \endcases\tag23
$$ 
For the 0'th and ($n+1$)'th face maps we further have $d_0\ho
l0(a_0\tens a_n)=\sigma_{lk}(a_0)\ot\Delta_N(a_1\tens 
a_n)$ and $d_{n+1}\ho ln(a_0\tens a_n)=\sigma_{lk}(a_0\tens a_n)$. Setting
 
$$
h^N_l=\sum_{i=0}^n(-1)^i\ho li\tag24
$$
gives the homotopy relations
$$
\sigma_{lk}(a_0\tens a_n)=\sigma_{lk}(a_0)\ot\Delta_N(a_1\tens a_n)+(\bd
h^N_l+h^N_l\bd)(a_0\tens a_n),\tag25
$$
see Lemma 1.0.9 of [L]. 

Let $\tau\in A_k^{\ot n}$ and suppose that
$\bd\tau=0$. Using Proposition 3.1 we have
$$
\sigma_{lk}(\tau)-\tau-\bd
h^N_l(\tau)=\bd(\int_0^1\alpha_t(s\tau)dt-h^N_l\tau),\tag26
$$
so 
$$
\aligned
\vert\phi(\sigma_{lk}(\tau)-\bd h^N_l(\tau))-\phi(\tau)\vert&=\vert
D(\int_0^1\alpha_t(s\tau)dt-h^N_l\tau)\vert\\
&\leq C\Vert D\Vert\Vert\tau\Vert.
\endaligned\tag27
$$
By the uniform bound on diagonal matrices given by condition M1 
the constant $C>0$ depends only on $n$. We next get
$$
\vert\phi(\frac1N\sum_{l=0}^{N-1}(\sigma_{lk}-\bd
h^N_l)\tau)-\phi(\tau)\vert\leq C\Vert D\Vert\Vert\tau\Vert.\tag28
$$
In particular
$$
(\phi(\frac1N\sum_{l=0}^{N-1}(\sigma_{lk}-\bd h^N_l)\tau))_{N\in\Bbb
N}\in\ell_\infty.\tag29
$$
Define
$$
\brf(\tau)=\LIM N(\phi(\frac1N\sum_{l=0}^{N-1}(\sigma_{lk}-\bd
h^N_l)\tau))\tag30
$$
and set $\psi_k=\phi-\brf$. Sofar we have shown that
$\Vert\psi_k(\tau)\Vert\leq C\Vert D\Vert\Vert \tau\Vert$ for $\bd
\tau=0$. In order to finish we must show that $\psi_k$
generates $D$ on $A_k$, that is we must show that $\brf(\bd\tau)=0$
when $\tau\in A_k^{\ot(n+1)}$.
To this end we shall use an other homotopy relation. On chains with
entries from $A_k$ we have from (25)
$$
\sum_{l=0}^{N-1}\sigma_{lk}=\Delta_N+\bd(\sum_{l=0}^{N-1}h^N_l)+(\sum_{l=0}^{N-1}h^N_l)\bd,\tag31
$$
where we used
$\sum_{l=0}^{N-1}\sigma_{lk}(a_0)\otimes\Delta_N(a_1\tens
a_n)=\Delta_N(a_0\tens a_n)$. Applying this to $\bd\tau$ we get
$$
\frac1N\sum_{l=0}^{N-1}(\sigma_{lk}-\bd
h^N_l)\bd\tau=\frac1N\Delta_N(\bd\tau)=\frac1N\bd(\Delta_N(\tau)),\tag32
$$
where the last equality follows from (21). We further get
$$
\aligned
\vert\phi(\frac1N\sum_{l=0}^{N-1}(\sigma_{lk}-\bd
h^N_l)\bd\tau)\vert&=\vert\phi(\frac1N\bd(\Delta_N(\tau)))\vert\\
&=\vert\frac1N D(\Delta_N(\tau))\vert\\
&\leq\frac1NC\Vert D\Vert\Vert\tau\Vert,
\endaligned\tag33
$$
where $C>0$ is the bound from M1. Letting $N\to\infty$ we get
$$
\brf(\bd\tau)=\LIM N\phi(\frac1N\sum_{l=0}^{N-1}(\sigma_{lk}-\bd
h^N_l)\bd\tau)=0,\tag34
$$
thus finishing the proof.
\enddemo
\remark{3.4 Remark} The use of H-unitality has in essence only involved
excission concerning extensions
$$
0@>>>\fA@>>>\fA\oplus\langle U\rangle@>>>\langle U\rangle@>>>0
$$
where $\langle U\rangle$ is a 2-dimensional Banach algebra with
$U^2=1$.
\endremark

\subheading{4. Some Applications}

First we note some well-known conditions to ensure partial
fulfillment of Condition M2.

\proclaim{4.1 Proposition} For an algebra $A$:
\flushpar\rm{(1)} If $A$ has a {\it splitting of the multiplication}, that is the
  multiplication map $\mu\colon A\ot A\to A$ has an  A-bimodule left
  inverse, then A is H-unital and $H_n(A,A)=H^n(A,A')=\{0\}$ for all $n\geq1$.
\flushpar\rm{(2)} If $A$ has local left or right units (i.e. for each finite set
$\text{\rm M}\subseteq A$ there is $e_{\text{\rm M}}\in A$ so that
$e_{\text{\rm M}}a=a,\;a\in \text{\rm M}$ or $ae_{\text{\rm
M}}=a,\;a\in \text{\rm M}$), then $A$ is H-unital.

For a Banach algebra $\fA$:
\flushpar\rm{(3)} If $\fA$ has a bounded one-sided approximate identity, then $\fA$ is strongly H-unital.
\endproclaim
\demo{Proof} (1): If $\fs\colon A\to A\ot A$ is a splitting of the multiplication, then $\fs\ot \bold1\colon C_{\bullet}(A)\to C_{\bullet+1}(A)$ is a contracting homotopy of the Hochschild-complex. The rest is contained in Proposition 2 and Proposition 5 of [W1].
\enddemo
Our first application is Banach algebras of compact operators on
vector-valued sequence spaces. Let $X$ be a Banach space and let $1\leq
 p<\infty$. For the Banach space of $p-$summable sequences,
$\ell_p(X)$, we consider the Banach algebra $\Cal A(\ell_p(X))$ of
appproximable operators.

\proclaim{4.2 Theorem} If $\A{\ell_p(X)},\;1\leq p<\infty$ is strongly
H-unital,  
then it is simplicially trivial.
\endproclaim
\demo{Proof} We show that $\A{\ell_p(X)}$ satisfies Conditions M1, M2,
and M3. Consider the shifts 
$$
\aligned &S(x_1,x_2,\dots)=(0,x_1,x_2,\dots),\\
&R(x_1,x_2,\dots)=(x_2,x_3,\dots)\qquad(x_1,x_2,\dots)\in\ell_p(X),
\endaligned\tag35
$$ 
and let $A=\Cal F(\ell_p(X))$, the algebra of finite rank operators on
$\ell_p(X)$. Then, by a well-known construction, $A_\infty$ as defined
in Condition M2, has a splitting of the multiplication: For $x\in X,
x^*\in X^*$ and $i,j\in\Bbb N$ consider the rank 1 operator given by
the elementary matrix $E_{ij}(x\ot x^*)$. Then
$A_\infty=\operatorname{span}\{E_{ij}(x\ot x^*)\mid x\in X,x^*\in X^*,
i,j\in\Bbb N\}$. Chose $e\in X,e^*\in X^*$ with $e^*(e)=1$. A
splitting of the multiplication is given by $E_{ij}(x\ot x^*)\mapsto
E_{i1}(x\ot e^*)\ot E_{1j}(e\ot x^*)$. By Proposition 4.1 and since
$S^kR^k@>>>0$ in the strong operator topology, the rest is clear.
\enddemo
For a Hilbert space $H$  it was proved in [W1] that  $\A H$, $\cB(H)$,
and the Calkin algebra 
$\Cal C(H)$ are all simplicially trivial. The same argument gives
\proclaim{4.3 Corollary} If $X$ has the bounded approximation property,
then  $\A{\ell_p(X)}$, $\cB(\ell_p(X))$, and the Calkin algebra $\Cal
C(\ell_p(X))$ are all simplicially trivial.
\endproclaim
\demo{Proof} When $X$ has the bounded approximation property, then so
has $\ell_p(X)$. It follows that $\A{\ell_p(X)}$ is a strongly
H-unital Banch algebra, and the sequence
$\nul@>>>\A{\lp X}@>>>\cB(\lp X)@>>>\Cal C(\lp X)@>>>\nul$ is weakly
split. It was noted in [D,G\&G] that $\cB(\lp X)$ is simplicially trivial. The
rest follows from Theorem 4.2 and excission.
\enddemo

If $(\mu,\Omega,\Sigma)$ is a measure space such that
$L_1(\mu,\Omega,\Sigma)$ is finite dimensional
(i.e\. $(\mu,\Omega,\Sigma)$ is {\it equivalent to a finite set}), then $L_p(\mu,X)$ is
isomorphic to a finite direct sum of copies of $X$ for all $p$. Hence,
if $\A X$ is strongly H-unital, Morita equivalence gives that
$\cH^{\bullet}(\A X,\A
X^*)\cong\cH^{\bullet}(\A{L_p(\mu,X)},\A{L_p(\mu,X)}^*)$. 

For infinite dimensional $L_p$-spaces we get

\proclaim{4.3 Corollary} Let $(\mu,\Omega,\Sigma)$ be a measure space
which is not equivalent to a finite set, and let $X$ be a Banach space
with the bounded approximation property. Then $\A{L_p(\mu,X)}$ is
simplicially trivial.
\endproclaim
\demo{Proof} In [G3] it is proved that $\A{\lp X}$ and $\A{L_p(\mu, X)}$
are Morita equivalent. For Banach algebras with bounded one-sided
approximate identities bounded Hoch\-schild homology with coefficients in
essential modules is Morita invariant [G1]. The result follows then from
Theorem  4.2.
\enddemo

\enddocument